\documentclass[pra,reprint]{revtex4-1}
\usepackage{graphicx}
\usepackage{appendix}
\usepackage{braket}
\usepackage{rotating}
\usepackage{amssymb}
\usepackage{amsmath}
\usepackage{txfonts}
\usepackage{bm}
\usepackage{color}
\usepackage{multirow}
\usepackage{booktabs}
\usepackage{dcolumn}
\usepackage{subcaption}

\begin{document}

\title{Graph-theoretical estimates of the diameters of the Rubik's Cube groups}

\author{So \surname{Hirata}}
\email{sohirata@illinois.edu}
\affiliation{Department of Chemistry, University of Illinois at Urbana-Champaign, Urbana, Illinois 61801, USA}

\date{\today}

\begin{abstract}
A strict lower bound for the diameter of a symmetric graph is proposed, which is calculable with the order $n$ and other local parameters of the graph 
such as the degree $k\,(\geq 3)$, even girth $g\,(\geq 4)$, and number of $g$-cycles traversing a vertex, which are easily determined 
by inspecting a small portion of the graph (unless the girth is large). It is applied to the symmetric Cayley graphs of some Rubik's Cube groups 
of various sizes and metrics, yielding slightly tighter lower bounds of the diameters than those for random $k$-regular graphs proposed by Bollob\'{a}s and de la Vega. 
They range from 60\% to 77\% of the correct diameters of large-$n$ graphs.
\end{abstract}

\maketitle 

\section{Introduction}

\begin{figure}
  \includegraphics[width=0.6\columnwidth]{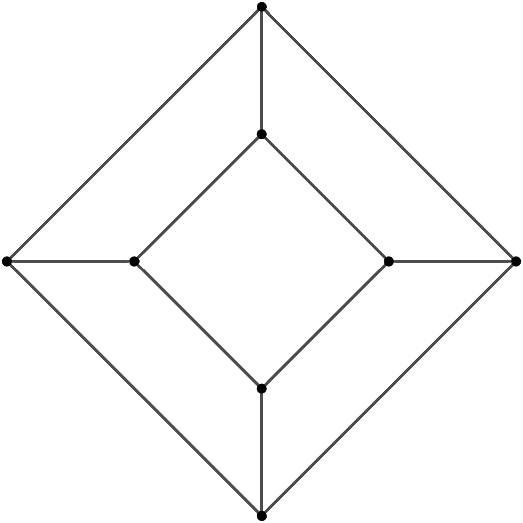}
\caption{Generalized Petersen graph $G(4,1)$ as the Cayley graph of the $3\times3\times3$ Rubik's Cube group in the square-slice-turn metric. Each vertex is a 
unique configuration (a pattern of face colors of the cubies while holding the orientation of the whole Cube fixed) reachable by a combination
of the $R^2L^{-2}$, $D^2U^{-2}$, and $B^2F^{-2}$ turns (in the Singmaster notation \cite{Singmasterbook}) from the solved configuration. 
Each edge connecting two vertexes is a turn (and its inverse) 
that brings one configuration denoted by one of the vertexes to the other.
The order of the graph is 8, the degree is 3, the girth ($g$) is 4, the diameter is 3, and the number of $g$-cycles traversing a vertex is 3.}
\label{fig:GPG4_1}
\end{figure}

\begin{figure}
  \includegraphics[width=0.6\columnwidth]{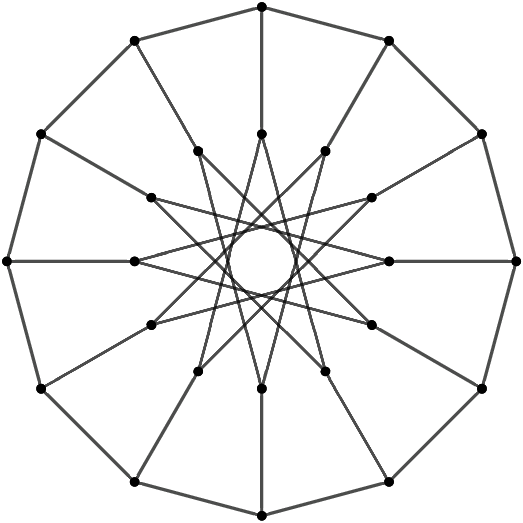}
\caption{Generalized Petersen graph $G(12,5)$ as the Cayley graph of the $2\times2\times2$ Rubik's Cube group in the square-turn metric.
The order is 24, the degree is 3, the girth ($g$) is 6, the diameter is 4, and the number of $g$-cycles per vertex is 3.}
\label{fig:GPG12_5}
\end{figure}

The fewest number of turns to solve the Rubik's Cube from its hardest initial configuration is colloquially referred to as God's Number \cite{Hofstadter,Singmasterbook,Joyner2008,Rokicki2013,Rokicki2014,cube20}.
It is equal to the diameter of the Cayley graph of the Rubik's Cube group. 
In the preceding article \cite{Hirata2024}, we estimated these diameters for the Cubes of various sizes in different metrics 
(a set of allowed combinations of rotations, each counted as one turn) using probabilistic logic. 
The estimates are accurate to within 2 of the actual values (when available) and are in exact agreement with the latter 
for both the 2$\times$2$\times$2 and 3$\times$3$\times$3 Cubes when the quarter-turn metric is adopted. 

The Cayley graphs of the Rubik's Cube group in some metrics are a symmetric graph \cite{Bollobasbook2,Biggsbook,Bollobasbook,Westbook}. For example,  the
3$\times$3$\times$3 Cube group in the square-slice-turn metric ($R^2L^{-2}$, $D^2U^{-2}$, and $B^2F^{-2}$ in the Singmaster 
notation \cite{Singmasterbook}) has the Cayley graph that is a cubic symmetric graph known as the generalized Petersen graph $G(4,1)$ shown in Fig.\ \ref{fig:GPG4_1} \cite{Westbook}.
The  Cayley graph of the 2$\times$2$\times$2 Cube group in the square-turn metric ($R^2$, $D^2$, and $B^2$) \cite{Singmasterbook} is 
another cubic symmetric graph, which is identified as the generalized Petersen graph $G(12,5)$ (Fig.\ \ref{fig:GPG12_5}) \cite{Westbook}.

The Cayley graphs of the 2$\times$2$\times$2 and 3$\times$3$\times$3 Cube groups in the quarter-turn metrics should
also be symmetric graphs because all configurations and all turns are symmetrically equivalent.
The Cayley graph of the 3$\times$3$\times$3 Cube group in the square metric is another symmetric graph, which turns out to be
locally isomorphic to the graph of the 2$\times$2$\times$2 Cube group in the quarter-turn metric, with the only difference being their orders (numbers of Cube configurations). 

In this article, we introduce the formula for a strict lower bound for the diameter of a general symmetric graph derived with algebraic-graph-theoretical arguments \cite{Biggsbook,Westbook}. It can be evaluated with the order and some easily 
determined local parameters of the graph, such as the degree ($k$), girth ($g$), and number of $g$-cycles traversing a vertex. 
We then apply this formula to the diameters of the Cube groups of different sizes and metrics and compare the lower bounds with the probabilistic estimates \cite{Hirata2024}
and actual values as well as with the lower and upper bounds of the diameter of a random $k$-regular graph proposed by Bollob\'{a}s and de la Vega \cite{BollobasVega1982,Bollobasbook}. See 
the preceding article \cite{Hirata2024} for previous studies and nomenclature of the Rubik's Cube groups.

\section{Symmetric Cayley graphs\label{sec:cayley}}

Following the notation of Biggs \cite{Biggsbook}, a distance-partition graph $\Gamma_i(v)$ for any vertex $v$ of a symmetric graph $\Gamma$ 
of order $n$, degree $k\, (\geq 3)$, and even girth $g\, (\geq 4)$ is defined by its vertex set,
\begin{eqnarray}
V\Gamma_i(v) = \left\{ u \in V\Gamma\, |\,  \partial(u,v) = i \right\},
\end{eqnarray}
where $\partial(u,v)$ is the distance between vertexes $u$ and $v$. Here, $V\Gamma_0(v) = \{ v\}$ and $|V\Gamma_0(v)| = 1$ as well as
\begin{eqnarray}
|V\Gamma_i(1)| = |V\Gamma_i(2)| = \dots = |V\Gamma_i(n)| \geq 1, 
\end{eqnarray}
at any $i$ ($0 \leq i \leq d$) because $\Gamma$ is a symmetric graph. The $V\Gamma_d(v)$ is the set of the antipodes of vertex $v$.

Its corresponding edge set is defined by
\begin{eqnarray}
E\Gamma_i(v) = \left\{ \{ u,v \} \in E\Gamma\, |\, u \in V\Gamma_{i}, v \in V\Gamma_{i-1} \right\},
\end{eqnarray}
with $E\Gamma_0(v) = \emptyset$ and $|E\Gamma_0(v)| = 0$ as well as
\begin{eqnarray}
|E\Gamma_i(1)| = |E\Gamma_i(2)| = \dots = |E\Gamma_i(n)| \geq 1,
\end{eqnarray}
at any $i$ ($1 \leq i \leq d$). If $\Gamma$ has a diameter $d$,  
\begin{eqnarray}
|E\Gamma_{d+1}(v)| = 0,
\end{eqnarray} 
for any $v$. Henceforth, we drop the spurious argument $v$.

We can easily show
\begin{eqnarray}
\sum_{i=0}^{d} |V\Gamma_i| &=& n, \label{prop1} \\
\sum_{i=1}^{d} |E\Gamma_i| &=& \frac{nk}{2}, \label{prop2}
\end{eqnarray}
as well as
\begin{eqnarray}
\sum_{i=0}^{d} (-1)^i |V\Gamma_i| &=& 0, \label{prop3}
\end{eqnarray}
using $|E\Gamma_{i}| + |E\Gamma_{i+1}| = k|V\Gamma_i|$ for the last identity. 

Let $\bm{d}$ be an array of orders of the distance-partition graph $V\Gamma_i$, i.e.,
\begin{eqnarray}
\bm{d} = \left\{ |V\Gamma_0|, |V\Gamma_1|, \dots, |V\Gamma_d| \right\} .
\end{eqnarray}
The $d_i$ ($0 \leq i \leq d$) is therefore the number of vertexes with distance $i$ from an arbitrary origin vertex; $d_0 = 1$ and $d_1 = k$. 
Let $\bm{r}$ be an array of {\it branching ratios} defined by
\begin{eqnarray}
\bm{r} = \left\{ \frac{|V\Gamma_1|}{|V\Gamma_0|}, \frac{|V\Gamma_2|}{|V\Gamma_1|}, \dots, \frac{|V\Gamma_d|}{|V\Gamma_{d-1}|} \right\} .
\end{eqnarray}
The $r_i$ ($1 \leq i \leq d$) is the rate of growth in the number of vertexes with increasing distance, i.e., $r_1 = k$, $r_2 \leq k-1$, and generally,
\begin{eqnarray}
r_{i+1} \leq r_{i}. \label{prop4} 
\end{eqnarray}

The above equations can be rationalized as follows: From the origin vertex (distance 0), $k$ edges emanate and reach $k$ new vertexes at distance 1 ($r_1 = k$). 
From each of the $k$ vertexes at distance 1, one edge connects the origin vertex, but the remaining $k-1$ edges spawn at most $k-1$ new vertexes at distance 2 ($r_2 \leq k-1$).
The equal sign applies when the girth is greater than 4. 
When the girth is equal to 4, some of the vertexes are a shared destination of two edges and therefore the number of new vertexes at distance 2 is less ($r_2 < k-1$).
In subsequent steps, possible existence of longer cycles decreases the branching ratio even further ($r_{i+1} \leq r_{i}$).

\section{Diameters of symmetric graphs}

In a symmetric graph with order $n$ and degree $k\, (\geq 3$), 
let $g\, (\geq 4$) be the even girth, $d$ be the diameter, and $\eta$ be the number of cycles of length $g$ traversing a vertex.
We seek a lower bound for $d$, whose precise determination is generally hard,  in terms of $n$ and local parameters of the graph---$k$, $g$, and $\eta$---which are usually 
readily available by insprection.  

The branching ratios $r$ up to distance $g/2$ are
\begin{eqnarray}
r_1 &=& k, \\
r_i &=& k-1 \,\,\, (2 \leq i \leq g/2-1),
\end{eqnarray}
and
\begin{eqnarray}
r_{g/2} &=& \frac{ k(k-1)^{g/2-1} - \eta }{k(k-1)^{g/2-2}} \equiv r_\text{max}. \label{growth}
\end{eqnarray}
As per Eq.\ (\ref{prop4}), this $r_{g/2}$ is the upper bound of $r$ for greater distances. We call it $r_\text{max}$, which can take any real positive value. 

An upper bound for $n$ is then obtained by assuming that the vertexes multiply at the same highest branching ratio of $r_\text{max}$ all the way to distance $d$.
\begin{eqnarray}
n_\text{max} &=& \left \lfloor 1+ k+k(k-1) + \dots + k(k-1)^{g/2-2} \right. \nonumber\\
&& \left. + k(k-1)^{g/2-2} r_\text{max} + \dots + k(k-1)^{g/2-2} r_\text{max}^{d-g/2+1} \right\rfloor. 
\end{eqnarray}

When $r_\text{max} \neq 1$, the above geometric series can be summed to
\begin{eqnarray}
n_\text{max}  &=& \left \lfloor n_0 + k(k-1)^{g/2-2} \frac{r_\text{max}^{d-g/2+2}-r_\text{max}}{r_\text{max}-1} \right\rfloor
\end{eqnarray}
with
\begin{eqnarray}
n_0 &=& \frac{k(k-1)^{g/2-1}-2}{k-2},
\end{eqnarray}
which is related to the Moore bound (and equal to the latter when $g=2d$ and $\eta=0$) \cite{Biggsbook,Bollobasbook,Westbook}.
When $r_\text{max}=1$, we instead have
\begin{eqnarray}
n_\text{max} &=& \left\lfloor n_0 + k(k-1)^{g/2-2} (d-g/2+1) \right\rfloor . 
\end{eqnarray}

The corresponding lower bound for $d$ is obtained by equating $n$ with $n_\text{max}$, which reads
\begin{eqnarray}
d_\text{min} &=&  \left\lceil \frac{g}{2} - 2 + \log_{r_\text{max}} \left\{ \frac{ (n-n_0)(r_\text{max}-1)}{k(k-1)^{g/2-2}} + r_\text{max} \right\} \right\rceil, \,\,\,(r_\text{max} \neq 1),  \label{dmin1}
\end{eqnarray}
or 
\begin{eqnarray}
d_\text{min} &=& \left\lceil \frac{g}{2} - 1 +  \frac{n-n_0}{k(k-1)^{g/2-2}} \right\rceil, \,\,\,(r_\text{max} = 1).\label{dmin2}
\end{eqnarray}

We have confirmed that they indeed give a lower bound for the diameter for each of the cubic symmetric graphs (with an even girth, conservatively assuming $\eta = 1$) listed in 
a Foster Census \cite{Conder2002}.

A similar argument can be used to determine an upper bound for $d$ by assuming that the number of vertexes grows at the slowest pace after distance $g/2$, i.e., 
just one vertex at each distance. Unfortunately, this gives a uselessly conservative upper bound and we shall not consider it any further.

In our previous study \cite{Hirata2024}, the diameter of a Rubik's Cube group is estimated by the following formula,
\begin{eqnarray}
d_\text{probab} &=& \log_{r_\text{max}} n + \frac{\log n}{r_\text{max}}. \label{probab}
\end{eqnarray}
which is neither a lower nor upper bound, but works for nonsymmetric Cayley graphs or for any turn metric that does not even form a mathematically well-defined group. 
It has been derived on the basis of a modified version of the coupon collector's problem of probability theory \cite{Hirata2024}.
Equation (\ref{probab}) can also be viewed as a probabilistic estimate of the diameter of a symmetric group, where $r_\text{max}$ is given by Eq.\ (\ref{growth}) (in Ref.\ \cite{Hirata2024}, 
a more accurate estimate of $r$ was computationally determined and used instead). The first term approximately counts the distance of the growth phase
of the distance array $\bm{d}$, where $r \approx r_\text{max}$
and the number of vertexes increases geometrically. The second term measures the distance of the subsequent rapid decay phase. Our graph-theoretical lower bound [Eq.\ (\ref{dmin1})] considers only the growth phase and its formula thus has a similar functional form with the first term of Eq.\ (\ref{probab}); they become nearly identical
in the $n \to \infty$ limit. 

Upper and lower bounds for the diameter of a random $k$-regular graph were obtained by Bollob\'{a}s and de la Vega \cite{BollobasVega1982,Bollobasbook}.
For almost every graph and a constant $\epsilon > 0$, the BV bounds are
\begin{eqnarray}
d_\text{min}^\text{BV} &=& \left\lfloor \log_{k-1} n \right\rfloor  
+ \left\lfloor \log_{k-1} \log n - \log_{k-1} \frac{6k}{k-2} \right\rfloor + 1,  \label{BVlower} \\
d_\text{max}^\text{BV} &=& \left\lceil {\log_{k-1} \left\{(2+\epsilon)\, k n \log n\right\} } \right\rceil + 1, \label{BVupper} 
\end{eqnarray}
asymptotically ($n \to \infty$). These are derived on the basis of probabilistic arguments not unlike Eq.\ (\ref{probab}). 
Hence, the lower bound [Eq.\ (\ref{BVlower})] bears resemblance with the first term of the probabilistic estimate [Eq.\ (\ref{probab})] or even $d_\text{min}$ [Eq.\ (\ref{dmin1})].
As $n \to \infty$, the first term of the former dominates and approaches the latter from below because $r_\text{max} \leq k-1$ [see Eq.\ (\ref{growth})]. 
In what follows, the BV bounds are computed with $\epsilon=0$. 

\section{Cayley graphs of Rubik's Cube groups}
 
 The reader is referred to the preceding article \cite{Hirata2024} for the terminologies, group/graph orders, and other details of the Rubik's Cube groups of various sizes and metrics. 
 Table \ref{table:summary} summarizes the results of applying above formulas. 
 
\begin{table*}
\caption{\label{table:summary} Symmetric Cayley graphs of the $n \times n \times n$ Rubik's Cube groups in various metrics.}
\begin{ruledtabular}
\begin{tabular}{cccddddddddd}
& & & & & & & \multicolumn{5}{c}{Diameter} \\ \cline{8-12}
\multicolumn{1}{c}{Cube} & \multicolumn{1}{c}{Metric} & \multicolumn{1}{c}{Order $n$}  & \multicolumn{1}{c}{Degree $k$}  & \multicolumn{1}{c}{Girth $g$}  & \multicolumn{1}{c}{$\eta$\footnotemark[1]} & \multicolumn{1}{c}{$r_\text{max}$\footnotemark[2]}& \multicolumn{1}{c}{$d$ (actual)\footnotemark[3]} & \multicolumn{1}{c}{$d_\text{min}$\footnotemark[4]} & \multicolumn{1}{c}{$d_\text{probab}$\footnotemark[5]} 
& \multicolumn{1}{c}{$d_\text{min}^\text{BV}$\footnotemark[6]}  & \multicolumn{1}{c}{$d_\text{max}^\text{BV}$\footnotemark[6]}  \\ \hline
3$\times$3$\times$3 & Square-slice & $8$ & 3 & 4 & 3 & 1.0 & 3 &  3 & \infty  &0 & 8 \\ 
2$\times$2$\times$2 & Square  & $24$ & 3 & 6 & 3 & 1.5 & 4 &  4 & 10.0 & 2 & 10  \\ 
2$\times$2$\times$2 & Quarter  & $3.67\times10^6$ & 6 & 4 & 3 &4.5 & 14 &  10 & 13.4 & 10 & 14 \\ 
3$\times$3$\times$3 & Sqaure & $6.63\times10^5$ & 6 & 4 & 3 & 4.5 & 15 &  9 &  11.9 & 9 & 13 \\ 
3$\times$3$\times$3 & Quarter & $4.33\times10^{19}$ & 12 & 4 & 18 & 9.5 & 26 &  20 & 24.8 & 19 & 23 \\ 
\end{tabular}
\end{ruledtabular} 
\footnotetext[1]{The number of $g$-cycles per vertex.}
\footnotetext[2]{Upper bound for the branching ratio at distance $g/2$ and greater [Eq.\ (\ref{growth})].}
\footnotetext[3]{The diameter of the 3$\times$3$\times$3 Cube group in the quarter-turn metric was determined computationally by Rokicki {\it et al.}\ \cite{cube20,Rokicki2013,Rokicki2014}.
Other diameters listed are easily determined computationally \cite{Hirata2024}.}
\footnotetext[4]{Graph-theoretical lower bound for the diameter $d$ [Eq. (\ref{dmin1}) or (\ref{dmin2})] (this work).}
\footnotetext[5]{Probabilistic estimate [Eq.\ (\ref{probab})] (the preceding article \cite{Hirata2024}).}
\footnotetext[6]{Asymptotic ($n \to \infty$) lower bound [Eq.\ (\ref{BVlower})] and upper bound ($\epsilon=0$) [Eq.\ (\ref{BVupper})]  of Bollob\'{a}s and de la Vega \cite{BollobasVega1982,Bollobasbook}.}
\end{table*}

\subsection{The 3$\times$3$\times$3 Cube in the square-slice-turn metric}
 
The Cayley graph of the 3$\times$3$\times$3 Rubik's Cube group in the square-slice-turn metric \cite{Singmasterbook} ($R^2L^{-2}$, $D^2U^{-2}$, and $B^2F^{-2}$)  is a cubic symmetric graph known as the generalized Petersen graph $G(4,1)$ depicted in Fig.\ \ref{fig:GPG4_1}.
The objective here is to determine the diameter $d$ of this graph from its global ($n$) and local parameters ($k$, $g$, and $\eta$)
that are easily determined without having the whole graph. (In this tiny example, we immediately know the answer:\ $d=3$.)

The order of the graph $n$ is a global parameter and is potentially hard to determine. However, for a Rubik's Cube group, it 
can be worked out by a careful combinatorial argument \cite{Hofstadter,Singmasterbook,Hirata2024} or with the aid of group theory \cite{Joyner2008}.
In the square-slice metric \cite{Singmasterbook}, each of the three center slices is allowed to independently take one of two orientations  and thus $n=2^3=8$. 
Other local parameters such as the degree $k$, girth $g$, and the number ($\eta$) of $g$-cycles per vertex can be gleaned by inspecting a 
small neighborhood of a vertex in the partial graph. From Fig.\ \ref{fig:GPG4_1}, we find $k=3$, $g=4$, and $\eta = 3$.

The distance array is $\bm{d} = \{ 1, 3, 3, 1\}$, which obeys Eqs.\ (\ref{prop1}) and (\ref{prop3}) as it must. 
The branching ratios are $\bm{r} = \{3,1,1/3\}$, satisfying Eq.\ (\ref{prop4}).

Using Eq.\ (\ref{growth}), $r_\text{max} = 1$, which corresponds to the predicted distance array of $\bm{d}=\{1, 3, 3, 3\}$. 
Equation (\ref{dmin2}) yields $d_\text{min} = \lceil 2.3 \rceil = 3$, which is a tight lower bound for the actual diameter of 3. 
The probabilistic estimate [Eq.\ (\ref{probab})] breaks down in this tiny case because $n=8$ is too few for probabilistic logic.
For the same reason, the BV lower and upper bounds [Eqs.\ (\ref{BVlower}) and (\ref{BVupper})] are too low and too high.

\subsection{The 2$\times$2$\times$2 Cube in the square-turn metric}

The Cayley graph of the 2$\times$2$\times$2 Rubik's Cube group in the square-turn metric \cite{Singmasterbook}  ($R^2$, $D^2$, and $B^2$) is also a cubic symmetric graph
known as the generalized Petersen graph $G(12,5)$ drawn in Fig.\ \ref{fig:GPG12_5}. It is also called the Nauru graph. 

In the absence of the whole Cayley graph, its order ($n$) can be determined  to be 24 by the following combinatorial argument: 
The {\it ulf} cubie \cite{Singmasterbook}  is held fixed as the orientation anchor. The {\it urf} cubie and three cubies ({\it ulf}, {\it drb}, {\it dlb}) that are diagonal from it can interchange their positions in
$4! = 24$ different ways. Their orientations and the other four cubies' positions and orientations are completely determined by the positions of the first four. Hence, $n=24$. 
Focusing on any one vertex and its small neighborhood of Fig.\ \ref{fig:GPG12_5}, we find $k=3$, $g=6$, and $\eta=3$. 

Substituting these values into Eq.\ (\ref{dmin1}), we obtain $d_\text{min} = \lceil 3.4 \rceil = 4$, which is a tight lower bound for the actual diameter of 4 (Table \ref{table:summary}). 
This graph is still too small for probabilistic arguments and $d_\text{probab} = d_\text{max}^\text{BV} = 10$ are too high, while $d_\text{min}^\text{BV} = 2$ is too low.

The distance array is $\bm{d} = \{ 1,3,6,9,5\}$, satisfying Eqs.\ (\ref{prop1}) and (\ref{prop3}), whereas the branching-ratio array is $\bm{r} = \{ 3,2,3/2,5/9\}$. The graph-theoretical estimations approximate them 
as $\bm{r} \approx \{ 3,2,3/2,3/2\}$ and thus $\bm{d} \approx \{ 1,3,6,9,13.5\}$ since $r_\text{max} = 1.5$. 
 
\subsection{The 2$\times$2$\times$2 Cube in the quarter-turn metric\label{sec:2Q}}

\begin{figure}
  \includegraphics[width=0.6\columnwidth]{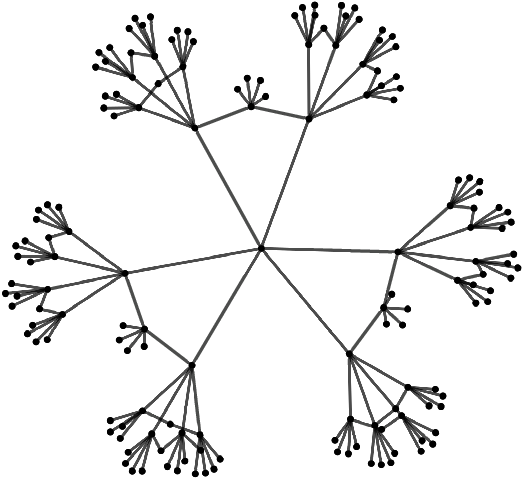}
\caption{A partial Cayley graph of the 2$\times$2$\times$2 Rubik's Cube group in the quarter-turn metric, which is symmetric with degree 6.
The order is $3.67 \times 10^6$, the girth ($g$) is 4, the diameter is 14, and the number of $g$-cycles per vertex is 3.
The vertexes on the periphery are further connected to other vertexes that are not shown.}
\label{fig:2x2x2Q}
\end{figure}

Figure \ref{fig:2x2x2Q} shows a local Cayley graph within distance 3 from an arbitrary origin vertex for the 2$\times$2$\times$2 Cube group in the quarter-turn metric \cite{Singmasterbook}. The order of the graph can be determined by combinatorial logic to be $n = 3.67 \times 10^6$ \cite{Hofstadter,Singmasterbook,Joyner2008,Hirata2024}.

Without drawing such a huge graph, based on the facts that all Cayley graphs are vertex-transitive and that 
all turns are symmetrically equivalent with one another, this Cayley graph should also be edge-transitive and thus symmetric \cite{Biggsbook,Westbook}.
This is supported by the appearance of a partial graph (within distance 3 from an arbitrarily chosen origin vertex) 
drawn in Fig.\ \ref{fig:2x2x2Q}, in which every vertex is equivalent and interchangeable with one another, and so is every edge.
This smallest portion of the Cayley graph is sufficient to determine $k=6$, $g=4$, and $\eta=3$. 

\begin{figure}[b]
  \includegraphics[width=\columnwidth]{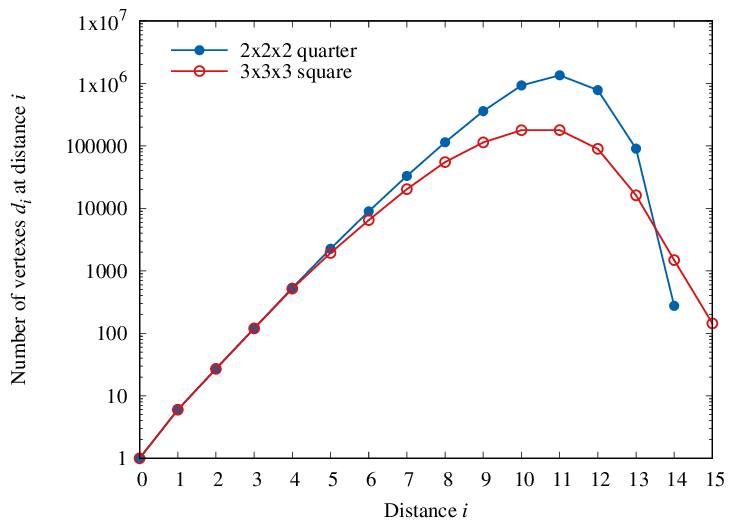}
\caption{The distance array $\bm{d}$ of the locally isomorphic Cayley graphs of the 2$\times$2$\times$2 Rubik's Cube group in the quarter-turn metric
and of the 3$\times$3$\times$3 Rubik's Cube group in the square-turn metric. The terminal point in each curve corresponds to the diameter.}
\label{fig:combined}
\end{figure}

Unlike these easily determined parameters, the diameter ($d$) is hard to come by. For the 2$\times$2$\times$2 Cube, 
a brute-force computation of the distance array $\bm{d}$ is possible \cite{Hirata2024,hirata_rubik}, as shown in Fig.\ \ref{fig:combined}. 
The computed $\bm{d}$ indicates that Eq.\ (\ref{prop3}) is obeyed and $d=14$. 
It is our objective, however, to estimate $d$ from the easily obtained local information of the graph as well as $n$.

Using Eq.\ (\ref{dmin1}), we find $d_\text{min} = \lceil 9.7 \rceil = 10$, which bounds the actual $d = 14$ from below as it should. 
The probabilistic estimate [Eq.\ (\ref{probab})] is $d_\text{probab} = 13.4$, which is much closer to the correct value of 14, as expected. 
The BV lower bound $d_\text{min}^\text{BV}=10$ [Eq.\ (\ref{BVlower})] is the same as $d_\text{min}$, which is not surprising considering the similarity in
their underlying derivations. The $d_\text{max}^\text{BV}=14$ is a tight upper bound and equal to the actual diameter ($d=14$) itself, but this is likely the result of error cancellation (see below). 
On the other hand, the good agreement between $d_\text{probab}=13.4$ and the actual $d$ is not accidental \cite{Hirata2024}.

\subsection{The 3$\times$3$\times$3 Cube in the square-turn metric}

\begin{figure}
  \includegraphics[width=0.6\columnwidth]{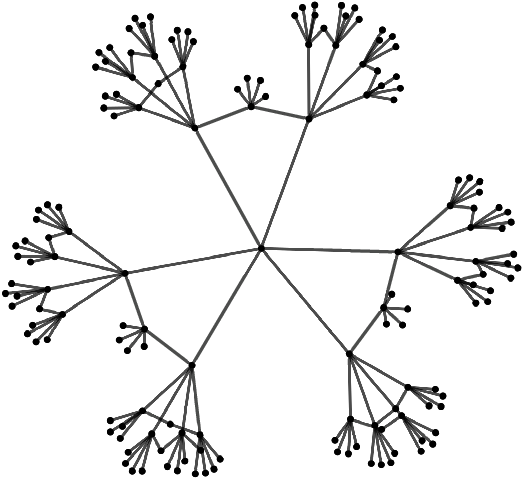}
\caption{A partial Cayley graph of the 3$\times$3$\times$3 Rubik's Cube group in the square-turn metric, which is symmetric with degree 6.
The order is $663\,552$, the girth ($g$) is 4, the diameter is 15, and the number of $g$-cycles per vertex is 3.}
\label{fig:3x3x3SQ}
\end{figure}

The square-turn metric of the 3$\times$3$\times$3 Rubik's Cube allows six 180$^\circ$ rotations denoted by $R^2$, $D^2$, $B^2$, $L^2$, $U^2$, and $F^2$
in the Singmaster notation \cite{Singmasterbook}. The order of the Cayley graph can be determined by group theory to be $663\,552$ \cite{Joyner2008,Hirata2024}.

What makes this group interesting is the fact that its Cayley graph is locally isomorphic with the graph of the 2$\times$2$\times$2 Cube group
in the quarter-turn metric (Sec.\ \ref{sec:2Q}). By ``locally isomorphic,'' we mean that their partial graphs within distance $g/2$ from an arbitrary origin vertex 
are pairwise isomorphic. 
The only difference is the orders of these graphs; one is 5.5 times smaller than the other.  
Figure \ref{fig:3x3x3SQ} shows the local (within distance 3) Cayley graph, which is indeed identical with Fig.\ \ref{fig:2x2x2Q}; it is 
a symmetric graph with $k=6$, $g=4$, and $\eta = 3$. Surprisingly, the diameter ($d=15$) of this graph with $n=6.63 \times 10^5$ is greater  
than the diameter ($d=14$) of the 2$\times$2$\times$2 Cube group in the quarter-turn metric with $n=3.67 \times 10^6$. 
With all the other parameters being the same, it is impossible to correctly predict this inversion by the simple graph-theoretical or probabilistic estimations presented
here and earlier \cite{Hirata2024,BollobasVega1982,Bollobasbook}.
These two graphs, therefore, underscore the diversity of symmetric graph structures differing only in their orders. 

Using Eq.\ (\ref{dmin1}), we obtain as the lower bound of this graph $d_\text{min}=9$, which is lower than $d_\text{min}=10$ of the 2$\times$2$\times$2 Cube group
in the quarter-turn metric. Likewise, the probabilistic estimate $d_\text{probab} = 11.9$ is also lower than that ($d_\text{probab} = 13.4$) of the 2$\times$2$\times$2 Cube group
in the quarter-turn metric and is uncharacteristically underestimated. The BV lower bound is 9 and equal to $d_\text{min}$ of this study. The BV upper bound of 13 fails in this case 
and is smaller than the actual diameter by two. This is likely because of our arbitrary choice of $\epsilon=0$ and  due to the fact that the asymptotic limit has not been reached with $n=663\,552$ (which is why
the exact agreement between $d_\text{max}^\text{BV}$ and $d$ in the 2$\times$2$\times$2 Cube in the quarter-turn metric is also likely accidental).
In all cases, the inversion is not reproduced. 

Figure \ref{fig:combined}  compares the distance arrays $\bm{d}$ of these two locally isomorphic Cayley 
graphs \cite{Hirata2024,hirata_rubik}. The two curves initially rise at the same rate for they have the same $r_\text{max}$. The curve for the 3$\times$3$\times$3 Cube group in
the square metric (the red open circles) begins to slow the growth earlier, which is expected because its $n$ is 5.5-times smaller. The peak in $\bm{d}$ occurs 
slightly earlier (by about one half distance unit) in the 3$\times$3$\times$3 Cube group in
the square metric with the smaller $n$. What causes the inversion is the decay phase; 
the $\bm{d}$ of the 3$\times$3$\times$3 Cube group in
the square metric decays unusually slowly (see other plots of $\bm{d}$ in Ref.\ \cite{Hirata2024} also), considerably 
delaying reaching the diameter. The cause of this slow decay is unknown.

\subsection{The 3$\times$3$\times$3 Cube in the quarter-turn metric}

\begin{figure}
  \includegraphics[width=0.6\columnwidth]{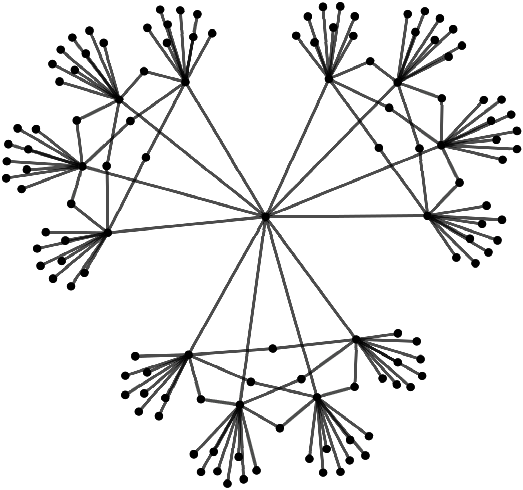}
\caption{A partial Cayley graph of the 3$\times$3$\times$3 Rubik's Cube group in the quarter-turn metric, which is symmetric with degree 12.
The order  is $4.33 \times 10^{19}$, the girth ($g$) is 4, the diameter is 26, and the number of $g$-cycles per vertex  is 18.}
\label{fig:3x3x3Q}
\end{figure}

Figure \ref{fig:3x3x3Q} is the local (within distance 2) Cayley graph of the 3$\times$3$\times$3 Rubik's Cube group in the quarter-turn metric \cite{Hirata2024}.
It is a symmetric graph with $k=12$, $g=4$, and $\eta = 18$. The order of the graph is $n=4.33\times10^{19}$ according to combinatorial logic \cite{Hofstadter,Singmasterbook,Joyner2008,Hirata2024}. The order of this group/graph is so high that an explicit enumeration of distance array $\bm{d}$ is nearly impossible. 
Nonetheless, by their Herculean computations, Rokicki {\it et al.}\  determined the diameter to be 26 \cite{Rokicki2013,Rokicki2014,cube20}.
 
Equation (\ref{dmin1}) places the lower bound for $d$ at $d_\text{min} = 20$, whereas the probabilistic estimation [Eq.\ (\ref{probab})] suggests $d_\text{probab} = 24.8$.
They are, respectively, consistent with and in reasonable agreement with the actual diameter of 26.  The $d_\text{min} = 20$ of this work is slightly tighter lower bound 
than the $d_\text{min}^\text{BV} = 19$, which is understandable because the former uses more information about the graph such as its symmetry and girth.
On the other hand, $d_\text{max}^\text{BV} = 23$ undershoots the diameter, indicating that the use of $\epsilon=0$ is inappropriate and/or $n=4.33 \times 10^{19}$ is still far from the asymptotic limit. 

\section{Conclusions}

Precisely determining the diameter of a large graph of the order $n$ in excess of millions or billions poses a major computational challenge. 
When $n$ can be known by a combinatorial or other argument and furthermore the graph has a high symmetry with its local parameters 
easily determined by inspection, one may estimate the diameter using an algebraic, geometrical, or probabilistic argument.
In the preceding article \cite{Hirata2024}, we proposed a probabilistic estimate of the diameter for a large Cayley graph, which reproduces the diameters of the Rubik's Cube groups
within 2 of the actual values or often exactly. Here, we introduce a strict lower bound for the diameter for a symmetric Cayley graph, calculable with $n$ and other easily
determined local parameters such as the degree $k$, girth $g$, and number $\eta$ of $g$-cycles per vertex. Similar, but more crude estimates of the diameters of the Rubik's Cube groups have been made \cite{Singmasterbook}, but this study is distinguished from them for presenting 
a rigorous lower bound also applicable to other symmetric graphs.
When applied to symmetric Cayley graphs of some 
Rubik's Cube groups, it yields tighter lower bounds than those for the diameter of a random $k$-regular graph proposed by
Bollob\'{a}s and de la Vega \cite{BollobasVega1982,Bollobasbook}. They are 60\% to 77\% of the correct values for large-$n$ graphs. 
 
\acknowledgments
The author is a Guggenheim Fellow of the John Simon Guggenheim Memorial Foundation. 

\bibliography{rubik.bib}

\begin{thebibliography}{14}%
\makeatletter
\providecommand \@ifxundefined [1]{%
 \@ifx{#1\undefined}
}%
\providecommand \@ifnum [1]{%
 \ifnum #1\expandafter \@firstoftwo
 \else \expandafter \@secondoftwo
 \fi
}%
\providecommand \@ifx [1]{%
 \ifx #1\expandafter \@firstoftwo
 \else \expandafter \@secondoftwo
 \fi
}%
\providecommand \natexlab [1]{#1}%
\providecommand \enquote  [1]{``#1''}%
\providecommand \bibnamefont  [1]{#1}%
\providecommand \bibfnamefont [1]{#1}%
\providecommand \citenamefont [1]{#1}%
\providecommand \href@noop [0]{\@secondoftwo}%
\providecommand \href [0]{\begingroup \@sanitize@url \@href}%
\providecommand \@href[1]{\@@startlink{#1}\@@href}%
\providecommand \@@href[1]{\endgroup#1\@@endlink}%
\providecommand \@sanitize@url [0]{\catcode `\\12\catcode `\$12\catcode
  `\&12\catcode `\#12\catcode `\^12\catcode `\_12\catcode `\%12\relax}%
\providecommand \@@startlink[1]{}%
\providecommand \@@endlink[0]{}%
\providecommand \url  [0]{\begingroup\@sanitize@url \@url }%
\providecommand \@url [1]{\endgroup\@href {#1}{\urlprefix }}%
\providecommand \urlprefix  [0]{URL }%
\providecommand \Eprint [0]{\href }%
\providecommand \doibase [0]{http://dx.doi.org/}%
\providecommand \selectlanguage [0]{\@gobble}%
\providecommand \bibinfo  [0]{\@secondoftwo}%
\providecommand \bibfield  [0]{\@secondoftwo}%
\providecommand \translation [1]{[#1]}%
\providecommand \BibitemOpen [0]{}%
\providecommand \bibitemStop [0]{}%
\providecommand \bibitemNoStop [0]{.\EOS\space}%
\providecommand \EOS [0]{\spacefactor3000\relax}%
\providecommand \BibitemShut  [1]{\csname bibitem#1\endcsname}%
\let\auto@bib@innerbib\@empty
\bibitem [{\citenamefont {Singmaster}(1981)}]{Singmasterbook}%
  \BibitemOpen
  \bibfield  {author} {\bibinfo {author} {\bibfnamefont {D.}~\bibnamefont
  {Singmaster}},\ }\href@noop {} {\emph {\bibinfo {title} {Note on {Rubik's}
  {Magic} {Cube}}}}\ (\bibinfo  {publisher} {Enslow Publishers},\ \bibinfo
  {address} {Hillside, NJ},\ \bibinfo {year} {1981})\BibitemShut {NoStop}%
\bibitem [{\citenamefont {Hofstadter}(1981)}]{Hofstadter}%
  \BibitemOpen
  \bibfield  {author} {\bibinfo {author} {\bibfnamefont {D.~R.}\ \bibnamefont
  {Hofstadter}},\ }\href@noop {} {\bibfield  {journal} {\bibinfo  {journal}
  {Sci. Am.}\ }\textbf {\bibinfo {volume} {244}},\ \bibinfo {pages} {20}
  (\bibinfo {year} {1981})}\BibitemShut {NoStop}%
\bibitem [{\citenamefont {Joyner}(2008)}]{Joyner2008}%
  \BibitemOpen
  \bibfield  {author} {\bibinfo {author} {\bibfnamefont {D.}~\bibnamefont
  {Joyner}},\ }\href@noop {} {\emph {\bibinfo {title} {Adventures in Group
  Theory: Rubik's Cube, Merlin's Machine, \& Other Mathematical Toys}}},\
  \bibinfo {edition} {2nd}\ ed.\ (\bibinfo  {publisher} {The Johns Hopkins
  University Press},\ \bibinfo {address} {Baltimore, MD},\ \bibinfo {year}
  {2008})\BibitemShut {NoStop}%
\bibitem [{\citenamefont {Rokicki}\ \emph {et~al.}(2013)\citenamefont
  {Rokicki}, \citenamefont {Kociemba}, \citenamefont {Davidson},\ and\
  \citenamefont {Dethridge}}]{Rokicki2013}%
  \BibitemOpen
  \bibfield  {author} {\bibinfo {author} {\bibfnamefont {T.}~\bibnamefont
  {Rokicki}}, \bibinfo {author} {\bibfnamefont {H.}~\bibnamefont {Kociemba}},
  \bibinfo {author} {\bibfnamefont {M.}~\bibnamefont {Davidson}}, \ and\
  \bibinfo {author} {\bibfnamefont {J.}~\bibnamefont {Dethridge}},\ }\href@noop
  {} {\bibfield  {journal} {\bibinfo  {journal} {{SIAM} J. Discrete Math.}\
  }\textbf {\bibinfo {volume} {27}},\ \bibinfo {pages} {1082} (\bibinfo {year}
  {2013})}\BibitemShut {NoStop}%
\bibitem [{\citenamefont {Rokicki}\ \emph {et~al.}(2014)\citenamefont
  {Rokicki}, \citenamefont {Kociemba}, \citenamefont {Davidson},\ and\
  \citenamefont {Dethridge}}]{Rokicki2014}%
  \BibitemOpen
  \bibfield  {author} {\bibinfo {author} {\bibfnamefont {T.}~\bibnamefont
  {Rokicki}}, \bibinfo {author} {\bibfnamefont {H.}~\bibnamefont {Kociemba}},
  \bibinfo {author} {\bibfnamefont {M.}~\bibnamefont {Davidson}}, \ and\
  \bibinfo {author} {\bibfnamefont {J.}~\bibnamefont {Dethridge}},\ }\href@noop
  {} {\bibfield  {journal} {\bibinfo  {journal} {{SIAM} Rev.}\ }\textbf
  {\bibinfo {volume} {56}},\ \bibinfo {pages} {645} (\bibinfo {year}
  {2014})}\BibitemShut {NoStop}%
\bibitem [{\citenamefont {Rokicki}\ \emph {et~al.}(2024)\citenamefont
  {Rokicki}, \citenamefont {Kociemba}, \citenamefont {Davidson},\ and\
  \citenamefont {Dethridge}}]{cube20}%
  \BibitemOpen
  \bibfield  {author} {\bibinfo {author} {\bibfnamefont {T.}~\bibnamefont
  {Rokicki}}, \bibinfo {author} {\bibfnamefont {H.}~\bibnamefont {Kociemba}},
  \bibinfo {author} {\bibfnamefont {M.}~\bibnamefont {Davidson}}, \ and\
  \bibinfo {author} {\bibfnamefont {J.}~\bibnamefont {Dethridge}},\ }\href@noop
  {} {\enquote {\bibinfo {title} {{God}'s {Number} is 20},}\ }\bibinfo
  {howpublished} {\url{https://www.cube20.org}} (\bibinfo {year} {accessed on
  08/27/2024})\BibitemShut {NoStop}%
\bibitem [{\citenamefont {Hirata}(2024{\natexlab{a}})}]{Hirata2024}%
  \BibitemOpen
  \bibfield  {author} {\bibinfo {author} {\bibfnamefont {S.}~\bibnamefont
  {Hirata}},\ }\href@noop {} {\enquote {\bibinfo {title} {Probabilistic
  estimates of the diameters of the {Rubik's Cube} groups},}\ }\bibinfo
  {howpublished} {arXiv:2404.07337} (\bibinfo {year}
  {2024}{\natexlab{a}})\BibitemShut {NoStop}%
\bibitem [{\citenamefont {Bollob\'{a}s}(1978)}]{Bollobasbook2}%
  \BibitemOpen
  \bibfield  {author} {\bibinfo {author} {\bibfnamefont {B.}~\bibnamefont
  {Bollob\'{a}s}},\ }\href@noop {} {\emph {\bibinfo {title} {Extremal Graph
  Theory}}}\ (\bibinfo  {publisher} {Academic Press},\ \bibinfo {address}
  {London},\ \bibinfo {year} {1978})\BibitemShut {NoStop}%
\bibitem [{\citenamefont {Biggs}(1993)}]{Biggsbook}%
  \BibitemOpen
  \bibfield  {author} {\bibinfo {author} {\bibfnamefont {N.}~\bibnamefont
  {Biggs}},\ }\href@noop {} {\emph {\bibinfo {title} {Algebraic Graph
  Theory}}},\ \bibinfo {edition} {2nd}\ ed.\ (\bibinfo  {publisher} {Cambridge
  University Press},\ \bibinfo {year} {1993})\BibitemShut {NoStop}%
\bibitem [{\citenamefont {Bollob\'{a}s}(2001)}]{Bollobasbook}%
  \BibitemOpen
  \bibfield  {author} {\bibinfo {author} {\bibfnamefont {B.}~\bibnamefont
  {Bollob\'{a}s}},\ }\href@noop {} {\emph {\bibinfo {title} {Random Graphs}}},\
  \bibinfo {edition} {2nd}\ ed.\ (\bibinfo  {publisher} {Cambridge University
  Press},\ \bibinfo {year} {2001})\BibitemShut {NoStop}%
\bibitem [{\citenamefont {West}(2019)}]{Westbook}%
  \BibitemOpen
  \bibfield  {author} {\bibinfo {author} {\bibfnamefont {D.~B.}\ \bibnamefont
  {West}},\ }\href@noop {} {\emph {\bibinfo {title} {Introduction to Graph
  Theory}}},\ \bibinfo {edition} {2nd}\ ed.\ (\bibinfo  {publisher} {Pearson},\
  \bibinfo {year} {2019})\BibitemShut {NoStop}%
\bibitem [{\citenamefont {Bollob\'{a}s}\ and\ \citenamefont {de~la
  Vega}(1982)}]{BollobasVega1982}%
  \BibitemOpen
  \bibfield  {author} {\bibinfo {author} {\bibfnamefont {B.}~\bibnamefont
  {Bollob\'{a}s}}\ and\ \bibinfo {author} {\bibfnamefont {W.~F.}\ \bibnamefont
  {de~la Vega}},\ }\href@noop {} {\bibfield  {journal} {\bibinfo  {journal}
  {Combinatorica}\ }\textbf {\bibinfo {volume} {2}},\ \bibinfo {pages} {125}
  (\bibinfo {year} {1982})}\BibitemShut {NoStop}%
\bibitem [{\citenamefont {Conder}\ and\ \citenamefont
  {Dobcs\'{a}nyi}(2002)}]{Conder2002}%
  \BibitemOpen
  \bibfield  {author} {\bibinfo {author} {\bibfnamefont {M.}~\bibnamefont
  {Conder}}\ and\ \bibinfo {author} {\bibfnamefont {P.}~\bibnamefont
  {Dobcs\'{a}nyi}},\ }\href@noop {} {\bibfield  {journal} {\bibinfo  {journal}
  {J. Comb. Math. Comb. Comput.}\ }\textbf {\bibinfo {volume} {40}},\ \bibinfo
  {pages} {41} (\bibinfo {year} {2002})}\BibitemShut {NoStop}%
\bibitem [{\citenamefont {Hirata}(2024{\natexlab{b}})}]{hirata_rubik}%
  \BibitemOpen
  \bibfield  {author} {\bibinfo {author} {\bibfnamefont {S.}~\bibnamefont
  {Hirata}},\ }\href@noop {} {\enquote {\bibinfo {title} {rubik},}\ }\bibinfo
  {howpublished} {\url{https://github.com/sohirata/rubik}} (\bibinfo {year}
  {accessed on 08/27/2024}{\natexlab{b}})\BibitemShut {NoStop}%
\end{thebibliography}%

\end{document}